\documentstyle{amsppt}
\pagewidth{6.0in}
\vsize8.0in
\parindent=6mm
\parskip=3pt
\baselineskip=16pt
\tolerance=10000
\hbadness=500
\NoRunningHeads

\topmatter

\title
Null form estimates for (1/2,1/2) symbols 
and local existence for a
quasilinear Dirichlet-wave equation 
\endtitle

\author
Hart F. Smith and Christopher D. Sogge
\endauthor

\thanks
Both authors are partially supported by the NSF.
\endthanks

\address
Department of Mathematics, University of Washington, Seattle, WA 98195
\endaddress

\email
hart\@math.washington.edu
\endemail

\address
Department of Mathematics, Johns Hopkins University, Baltimore, MD 21218
\endaddress

\email
sogge\@jhu.edu
\endemail



\abstract
We establish certain null form estimates of Klainerman-Machedon for
parametrices of variable coefficient wave equations for the convex obstacle
problem, and for wave equations with metrics of bounded curvature.  These
are then used to prove a local existence theorem for nonlinear
Dirichlet-wave equations outside of convex obstacles.

\endabstract

\endtopmatter

\document

\define\R{\Bbb R}

\head
1. Introduction
\endhead

\noindent
The purpose of this paper is to establish the following null form
estimate
$$
\bigl\|Q(du,dv)\bigr\|_{H^1(\Bbb R^{1+3}_{t,x})}\le C\,
\bigl(\,\|u_0\|_{H^2(\Bbb R^3)}+\|u_1\|_{H^1(\Bbb R^3)}\,\bigr)
\bigl(\,\|v_0\|_{H^2(\Bbb R^3)}+\|v_1\|_{H^1(\Bbb R^3)}\,\bigr)\,,
\tag1.1$$
for solutions $u$ and $v$ to the Cauchy problem for certain wave equations
$$
\cases
\partial^2_t u(t,x)=\Delta_{\bold g} u(t,x)\,,
\\
u(0,x)=u_0(x), \quad
\partial_t u(0,x)=u_1(x)\,.
\endcases
$$
The null form $Q$ may be of any one of the following forms
$$
\matrix\format\r&\;\;\;\c\;\;\;&\l\\
Q_0(du,dv)&=&\displaystyle\partial_tu(t,x)\,\partial_tv(t,x)-\sum_{i,j=1}^3
\bold g^{ij}(x)\partial_{x_i}u(t,x)\,\partial_{x_j}v(t,x)\,,\\
\\
Q_{\alpha\beta}(du,dv)&=&
\partial_{x_\alpha}u(t,x)\,\partial_{x_\beta}v(t,x)-
\partial_{x_\beta}u(t,x)\,\partial_{x_\alpha}v(t,x)\,,
\endmatrix
$$
where $x_\alpha$ and $x_\beta$ may represent $t$ or any $x_i\,.$
Here $\sum_{1\le i,j \le 3}{\bold g}^{ij}(x)\,d\xi_i\,d\xi_j$
denotes the cometric associated
with $\Delta_{\bold g}$.

For the Euclidean metric on $\R^3$, the estimate $(1.1)$ was established
globally by Klainerman and Machedon \cite{2}. For smooth
variable coefficient hyperbolic operators, 
local versions of (1.1) were
established by the second author in \cite{11}.

This paper is concerned with two new cases.
The first is the case that the wave equation is satisfied by
$u$ and $v$ for $x$ belonging to an open subset $\Omega\subset\R^3$
which has smooth boundary $\partial\Omega$, such that
$\partial\Omega\subset \R^3$ is
strictly geodesically concave with respect to $\bold g$.
We then assume that $u$ and $v$ satisfy
Dirichlet conditions on $\partial\Omega$,
$$
u(t,x)\bigr|_{x\in\partial\Omega}=0\,,\qquad
v(t,x)\bigr|_{x\in\partial\Omega}=0\,.
$$
In this case we prove $(1.1)$ for $t$ in a small time interval
and $x$ in the intersection of a small ball with $\Omega.$
We point out that when $\Omega$ is the complement of a strictly convex
obstacle in $\R^3$, with $\bold g$ the Euclidean metric, a partition
of unity argument, together with the global Euclidean estimates
of \cite{2}, implies $(1.1)$ globally in $x$,
for $t$ in any bounded interval.

The second case that our results apply to is where
$\bold g$ is a metric on a ball in $\R^3$, such that
the components of the Riemann curvature tensor of $\bold g$
are bounded measurable functions, and such that
the coordinate functions $x_i$ are harmonic with respect to
$\Delta_{\bold g}$. In such coordinates the metric
coefficients $\bold g_{ij}$ have second derivatives
belonging to $\text{BMO}(\R^3)$, and the geodesic flow
is uniquely determined and bilipschitz.
The solution operator for the wave equation in this situation
is studied in \cite{8}, \cite{9}.
It can be written as the composition of an operator of Fourier
integral type described below, with an operator which preserves
the Sobolev spaces $H^j(\R^3)\,,\; j=1,2\,.$ It then suffices
to establish mapping properties for the Fourier integral part,
which is the purpose of this paper.
The results of this paper will then imply that (1.1)
holds for such metrics provided that the norm is taken over
a set of unit size.

In both of the above cases,
the problem is reduced to establishing the following estimate
$$
\bigl\|Q(d\,T\!f,d\,T\!g)\bigr\|_{L^2(\Bbb R^{1+3}_{t,x})}\le
C\,\|f\|_{H^1(\Bbb R^3)}\|g\|_{H^2(\Bbb R^3)}\,,
\tag1.2$$
for an appropriate parametrix $T$ of order $0$.

For the obstacle problem, the main part of the parametrix takes the form
$$
T\!f(t,x)=\sum_{\pm}\int e^{i\varphi^\pm(t,x,\xi)}\,a^\pm(t,x,\xi)\,
\widehat{f}(\xi)\,d\xi\,,
$$
where the phases $\varphi^\pm(t,x,\xi)$ satisfies the eikonal equation
$$
\matrix\format\r\;\;&\c&\;\;\l\\
|\partial_t\varphi^\pm(t,x,\xi)|&=&\pm\|d_x\varphi(t,x,\xi)\|_{\bold g}\,,\\
\\
\varphi^\pm(0,x,\xi)&=&\langle x,\xi\rangle\,,\\
\endmatrix
\tag1.3$$
and the symbols, which vanish for $|\xi|\le 1$,
satisfy the following
modified $S^0_{\frac 23,\frac 13}$ estimates
$$
\bigl|\langle\xi,\partial_\xi\rangle^N\partial_{t,x}^\beta
\partial_\xi^\alpha a^\pm(t,x,\xi)\bigr|\le C_{N,\alpha,\beta}
\,\bigl(1+|\xi|\bigr)^{\frac{|\beta|}3-\frac{2|\alpha|}3}\,.
\tag1.4$$
There is also a ``diffractive'' term,
the estimation of which requires a modification of the argument for
the main term, as will be discussed in Section 4.

In the case of the wave equation for metrics of bounded curvature tensor,
the parametrix is more complicated. It takes the form
$$
T\!f(t,x)=\sum_{\pm}\sum_{k=1}^\infty
\int e^{i\varphi_k^\pm(t,x,\xi)}\,a^\pm_k(t,x,\xi)\,\widehat{f_k}
(\xi)\,d\xi\,,
\tag1.5
$$
where $\widehat{f}(\xi)=\sum_{k=0}^\infty\widehat{f_k}(\xi)\,,$
and for $k\ge 1$
the support of $\widehat{f_k}(\xi)$ lies in $2^{k-1}\le |\xi|\le 2^{k+1}\,.$

The phases $\varphi^\pm_k$, 
each of which is homogeneous of degree 1 in $\xi$,
satisfy the eikonal equation $(1.3)$ for a corresponding family
of metrics $\bold g_k$, where $\bold g_k$ is a sequence of smooth metrics
approximating the singular metric $\bold g$. This sequence of metrics
satisfies the estimates
$$
\bigl|\partial_x^\alpha \bold g_k(x)\bigr|\le
\cases
C\,,& |\alpha|\le 1\,,\\
C\,k\,,& |\alpha|=2\,,\\
C_\alpha\,2^{k(|\alpha|-2)/2}\,,&|\alpha|\ge 3\,.
\endcases
\tag1.6$$
It also satisfies
$$
\matrix\format\r\;\;&\c&\;\;\l\\
\bigl|\,\bold g_k(x)- \bold g(x)\bigr| & \le & C\,2^{-k}\,,\\
\\
\bigl|\,\nabla_x \bold g_k(x)-\nabla_x \bold g(x)\bigr| & \le &
C\,2^{-\frac k2}\,.
\endmatrix
\tag1.7$$
The sequence of phases satisfies corresponding estimates
$$
\sup_{|\xi|=1}\,\bigl|\partial_{t,x,\xi}^\alpha \varphi_k(t,x,\xi)\bigr|\le
\cases C\,,& |\alpha|\le 2\,,\\C\,2^{k(|\alpha|-2)/2}\,,&|\alpha|\ge 2\,.
\endcases
\tag1.8$$
It also satisfies, for $k\ge j$,
$$
\matrix\format\r\;\;&\c&\;\;\l\\
\displaystyle
\sup_{|\xi|=1}\,\bigl|\,\varphi^\pm_k(t,x,\xi)-\varphi^\pm_j(t,x,\xi)\bigr|
& \le & C\,2^{-j}\,,\\
\\
\displaystyle
\sup_{|\xi|=1}\,\bigl|\,\nabla_{t,x,\xi}\varphi^\pm_k(t,x,\xi)-
\nabla_{t,x,\xi}\varphi^\pm_j(t,x,\xi)\bigr| & \le & C\,2^{-\frac j2}\,.
\endmatrix
\tag1.9$$
Finally, the symbols satisfy the following
modified $S^0_{\frac 12,\frac 12}$ estimates,
$$
\bigl|\langle\xi,\partial_\xi\rangle^N\partial_{t,x}^\beta\partial_\xi^\alpha
a^\pm_k(t,x,\xi)\bigr|\le C_{N,\alpha,\beta}
\,2^{k\left(\frac{|\beta|}2-\frac{|\alpha|}2\right)}\,.
\tag1.10$$

One of the main motivations for establishing the estimate 
(1.1) is that it gives local
existence results for nonlinear wave equations with null form nonlinearities.
Consider, for example, an $N$ component system of the form
$$\cases
\partial^2_tu-\Delta_{\bold g}u=F(u,du),\,\quad x\in\Omega\,,
\\
u(0,\cdot\,)=u_0, \quad \partial_tu(0,\cdot\,)=u_1,
\\
u(t,\cdot\,)|_{\partial\Omega}=0,
\endcases
\tag1.11
$$
where $\Omega$ has geodesically concave boundary as discussed above.
We assume that $F(u,du)=(F^1(u,du),\dots,F^N(u,du))$, and
$$F^i(u,du)=\sum_{j,k}a^i_{j,k}(t,x)\Gamma^i_{j,k}(u)B^i_{j,k}(du^j,du^k),$$
with $B^i_{j,k}$ being a null form associated with $\bold g$, $a^i_{j,k}\in
C^\infty(\Bbb R\times \Omega)$, and $\Gamma^i_{j,k}\in C^\infty(\Bbb C^N)$.

If $u$ is a solution of (1.11), 
then the vanishing of $u$ and $\partial_t u$ on $\partial\Omega$
imposes the following compatability conditions on the data,
$$
u_0(x)=u_1(x)=0, \quad \text{if } \, x\in \partial \Omega\,.
\tag1.12$$
Conversely, under the hypotheses (1.12) on the data, 
we shall be able to obtain the following local existence result,
generalizing results from \cite{2} and \cite{11}.

\proclaim{Theorem 1.1}  Suppose that $u_j\in H^{2-j}(\Omega)$, $j=0,1$, have
compact support and satisfy $(1.12)$.  Then there is a $T_*>0$ and a unique
solution $u\in H^2([0,T_*]\times \Omega)$ of $(1.11)$ verifying
$$\|Q(du^j,du^k)\|_{H^1([0,T_*]\times \Omega)}<\infty, \quad 1\le j,k \le
N.$$
\endproclaim

We will return to this theorem in section 4, in which we
also discuss the reduction of the estimate (1.1) for
the obstacle problem to that of (1.2), and handle
the diffractive term. The main work of this
paper, which occupies sections 2 and 3, is to establish
estimate (1.2) for parametrices of the above types. 
Since the main part of the parametrix associated with the
convex obstacle problem is a special case of
the type (1.5) that arises from bounded curvature metrics,
we shall consider parametrices of the type (1.5) in 
sections 2 and 3.

\head
2. Further Reductions
\endhead

\noindent
We begin by reducing the proof of estimate $(1.2)$ to consideration
of the case that $\widehat{f}(\xi)$
is supported in a dyadic annulus at scale $2^k$, and $\widehat{g}(\xi)$
is supported in a ball of radius $c\,2^k$, where one may choose
$c$ arbitrarily small but fixed. To do this,
we fix $\beta\in C^\infty_0((1/2,2))$ so that
$\sum_{-\infty}^\infty \beta(2^js)=1$, $s>0$.  We then set
$$\hat f_k(\xi)=\beta(|\xi|/2^k)\widehat f(\xi)$$
so that $f=\sum f_k$ and $\text{supp }\widehat f_k\subset
\{\xi: \, 2^{k-1}\le |\xi|\le 2^{k+1}\}$.  We then write
$$
\tilde f_j=\sum_{k< j+N}\,f_k\,,\qquad
\tilde g_k=\sum_{j\le k-N}\,g_j\,.
$$
where $N$ is a fixed number that is to be specified later.
Recalling that the symbol of $T$ vanishes for small $|\xi|$, we have the
following identity,
$$
Q(dT\!f,dT\!g)=
\sum_{j=0}^\infty Q(dT\!\tilde f_j,dT\!g_j)
+
\sum_{k=0}^\infty
Q(dT\!f_k,dT\tilde g_k)
=I+II\,.
$$
We consider $I$ first.  By
the Strichartz estimates, which hold for the parametrix $T$ by
\cite{10} and \cite{7}, \cite{8}, we may bound
$$
\matrix\format\r\;\;&\c&\;\;\l\\
\displaystyle
\sum_{j=0}^\infty\,
\bigl\|Q(d\,T\!\tilde f_j,d\,T\!g_j)\bigr\|_{L^2(\R^{1+3}_{t,x})}
&\le&
\displaystyle
C\,\sum_{j=0}^\infty \|\tilde f_j\|_{H^{\frac 32}(\R^3)}\,
\|g_j\|_{H^{\frac 32}(\R^3)}\\
\\
&\le&
\displaystyle
C\,\Bigl(\sum_{j=0}^\infty
2^{-j}\,\|\tilde f_j\|^2_{H^{\frac 32}(\R^3)}
\Bigr)^{\frac 12}
\,\|g\|_{H^2(\R^3)}\\
\\
&\le&
C\;\|f\|_{H^1(\R^3)}\,\|g\|_{H^2(\R^3)}\,.
\endmatrix
$$

It thus remains to estimate $II$.
To estimate its $L^2$ norm, we first observe that,
for $N$ large enough, the terms are essentially
mutually orthogonal over $k$. This follows by a simple integration
by parts argument, which yields
$$
\left|\int
Q(d\,T\!f_k,d\,T\tilde g_k)\,\overline{Q(d\,T\!f_{k'},d\,T\tilde
g_{k'})}\,dt\,dx\,
\right|
\le
C\,2^{-|k-k'|}\,\|f_k\|_2\,\|g_k\|_2\,\|f_{k'}\|_2\,\|g_{k'}\|_2\,,
$$
provided $|k-k'|\ge 3$.
Consequently,
$$
\Bigl\|\;\sum_k Q(d\,T\!f_k,d\,T\tilde g_k)\,\Bigr\|^2_{L^2(\R^{1+3}_{t,x})}
\le C
\sum_k \bigl\|\,Q(d\,T\!f_k,d\,T\tilde g_k)\bigr\|^2_{L^2(\R^{1+3}_{t,x})}
+C\,
\|f\|_{L^2(\R^3)}\,\|g\|_{L^2(\R^3)}\,.
$$
Thus, to establish $(1.2)$, it suffices to establish
the following estimate, uniformly
over $k$:
$$
\bigl\|\,Q(d\,T\!f_k,d\,T\tilde g_k)\bigr\|_{L^2(\R^{1+3}_{t,x})}
\le
C\,\|f_k\|_{H^1(\R^3)}\,\|\tilde g_k\|_{H^2(\R^3)}\,.
$$
Finally, by the first estimate in $(1.7)$, another application of the
Strichartz estimates shows that we may replace
the metric $\bold g$ in the form $Q_0$ by the metric $\bold g_k$.

By writing $T=T^++T^-$, there are essentially two terms to consider:
$Q(d\,T^+\!f,d\,T^+\!g)$, and $Q(d\,T^+\!f,d\,T^-\!g)\,.$
In what follows we consider the term $Q(d\,T^+\!f,d\,T^+\!g)$;
the arguments hold with minor modification for the latter
term. To simplify notation we use $\varphi_k(t,x,\xi)$ to denote
$\varphi^+_k(t,x,\xi)\,.$

In the formula for the operator $d\,T^+$, the terms where $d$ hits the
symbol $a(t,x,\xi)$ are easily handled by the Strichartz and
energy estimates; thus
it suffices to restrict attention to the term where $d$ falls on the phase.
Let
$$
q_{kj}(t,x,\xi,\eta)=
Q\bigl(d\varphi_k(t,x,\xi/|\xi|),d\varphi_j(t,x,\eta/|\eta|)\bigr)\,.
$$

The next reduction is to introduce polar coordinates for the
$\eta$ variable, $\eta=\rho\omega$, where $\rho\in\R^+$, and
$\omega\in S^2$, the unit two-sphere.
We now fix $k$ and introduce the operator $T^\omega=T^\omega_k$ given by
$$
T^\omega(f,g)=\sum_{j\le k-N}
\int e^{i\varphi_k(t,x,\xi)+i\rho\varphi_j(t,x,\omega)}
a_k(t,x,\xi)\,a_j(t,x,\rho\omega)\,q_{kj}(t,x,\xi,\omega)
\,\widehat{f}_k(\xi)\,\widehat{g}_j(\rho)\,d\xi\,d\rho\,,
$$
where $f\in L^2(\Bbb R^3)$ and $g\in L^2(\Bbb R)$.
A simple argument (see, e.g. \cite{11}),
now reduces the proof of $(1.2)$ to showing that, for
$g_\omega\in L^2(\Bbb R\times S^2)$, the following holds
$$
\Bigl\|\int T^\omega(f,g_\omega)\,d\omega\Bigr\|_{L^2(dx\,dt)}
\le C\,\|f\|_{L^2(\Bbb R^3)}\,\|g_\omega\|_{L^2(\Bbb R\times S^2)}\,,
\tag2.1$$
where the Fourier transforms of $f$ and $g_\omega$ are restricted as
above.

The next step, following \cite{1}, is to decompose phase space into regions
on which the null form symbol $q_{kj}$ is essentially constant.
Since the phases $\varphi_j$ depend on the scale $j$, this
cannot be expressed simply in terms of the angle of $\xi$ to $\eta$.
To proceed, we set
$$
\delta(l)=2^{l-\frac k4}\,,
$$
and if $\beta$ is as above, we write
$$
q^l_{kj}(t,x,\xi,\eta)
=\beta\Bigl(\,\delta(l)^{-1}\times\text{angle}
\bigl[d_x\varphi_k(t,x,\xi),d_x\varphi_j(t,x,\eta)\bigr]\,\Bigr)\,
q_{kj}(t,x,\xi,\eta)\,.
$$
We then have
$$
q_{kj}(t,x,\xi,\eta)=q^0_{kj}(t,x,\xi,\eta)+
\sum_{l=1}^\infty\,q^l_{kj}(t,x,\xi,\eta)\,,
$$
where $q^0_{kj}(t,x,\xi,\eta)$ is supported in the region on which
the angle is bounded by $2\cdot 2^{-\frac k4}\,.$

Using this decomposition, we write
$T^\omega=\sum_l T^{l,\omega}\,,$
where (recall that $k$ is fixed)
$$\multline
T^{l,\omega}(f,g)=\\
\sum_{j\le k-N}
\int e^{i\varphi_k(t,x,\xi)+i\rho\varphi_j(t,x,\omega)}
a_k(t,x,\xi)\,a_j(t,x,\rho\omega)\,q^l_{kj}(t,x,\xi,\omega)
\,\widehat{f}_k(\xi)\,\widehat{g}_j(\rho)\,
\,d\xi\,d\rho\,.
\endmultline
\tag2.2$$
By $(1.8)$, and $(1.6)$ in the case of the null form $Q_0$
(recall that the metric $\bold g$ is replaced by $\bold g_k$),
the following estimates are valid for $j\le k$:
$$
\matrix\format\r\;\;&\c&\;\;\l\\
\bigl|\langle\xi,\partial_\xi\rangle^N\partial_{t,x}^\beta\partial_\xi^\alpha
q^l_{kj}(t,x,\xi,\omega)\bigr| & \le &
C_{N,\alpha,\beta}\,\delta(l)\,2^{\frac k2(|\beta|-|\alpha|)}\,,\\
\\
\bigl|\langle\xi,\partial_\xi\rangle^N\partial_{t,x}^\beta\partial_\xi^\alpha
\bigl(q^l_{kj}(t,x,\xi,\omega)-q^l_{kk}(t,x,\xi,\omega)\bigr)\bigr| & \le &
C_{N,\alpha,\beta}\,2^{-\frac j2}\,2^{\frac k2(|\beta|-|\alpha|)}\,.
\endmatrix
\tag2.3$$

For the next step, if $l$ is fixed, choose unit vectors 
$\xi^\mu \in S^2$ so
that the balls $B(\xi^\mu,\delta(l))$ cover $S^2$
with bounded overlap (independent of $\delta(l)$).
We then fix an associated partition of unity
$$
1=\sum_\mu \Psi^\mu(\xi), \qquad \xi\ne0
$$
consisting of $C^\infty (\Bbb R^3\backslash 0)$ functions that are
homogeneous of degree zero, and which satisfy
$$\text{supp }\Psi^\mu\cap S^2\subset B(\xi^\mu, 2\delta(l))\,,
\quad
D^\alpha \Psi^\mu(\xi)=O(\delta(l)^{-|\alpha|}) \quad\text{if}\quad
|\xi|=1\,.$$
If we then write
$$\widehat{f}(\xi)=\sum_\mu\widehat{f}_\mu(\xi)$$
we have the following

\proclaim{Lemma 2.1}
For fixed $N$ as above sufficiently large, the following holds
for $l\ge 0$,
$$
\Bigl\|\,\sum_\mu \int T^{l,\omega}(f_\mu,g_\omega)\,d\omega\,
\Bigr\|^2_{L^2(dx\,dt)}\le C\,
\sum_\mu \Bigl\|\, \int T^{l,\omega}(f_\mu,g_\omega)\,d\omega\,
\Bigr\|^2_{L^2(dx\,dt)}
+
C\,\|f\|_{L^2(\Bbb R^3)}\,\|g_\omega\|_{L^2(\Bbb R\times S^2)}\,.
$$
\endproclaim
\demo{Proof}
We shall show that if $C$ is a large constant and $|\xi^\mu-\xi^{\mu'}|\ge
C\,\delta(l)\,,$ then for any $M>0$,
$$
\int T^{l,\omega}(f_\mu,g)(t,x)\,
\overline{T^{l,\omega'}(f_{\mu'},g')(t,x)}\,dt\,dx\le
C_M\,2^{-kM}\,\|f_\mu\|_2\,\|f_{\mu'}\|_2\,\|g\|_2\,\|g'\|_2\,.
\tag2.4$$
This follows by considering the operator $(T^{l,\omega})^*\,
T^{l,\omega'}\,.$
If
$$q^l_{kj}(t,x,\xi,\rho\omega)\,q^l_{kj'}(t,x,\xi',\rho'\omega')
\,\widehat f_\mu(\xi)\,\widehat f_{\mu'}(\xi')\ne 0\,,$$
then the angle of
$\nabla_x\varphi_k(t,x,\xi)+\rho\nabla_x\varphi_j(t,x,\omega)$
to
$\nabla_x\varphi_k(t,x,\xi')+\rho'\nabla_x\varphi_{j'}(t,x,\omega')$
is bounded below by $\delta(l)$, provided $|\xi^\mu-\xi^{\mu'}|\ge
C\,\delta(l)$
for some large $C$,  $|\xi|, \, |\xi'|\approx 2^k$, and $\rho,\rho'\le
2^{k-N}$ with $N$ sufficiently large. On account of this,
$$
\bigl|
\bigl(\nabla_x\varphi_k(t,x,\xi)+\rho\nabla_x\varphi_j(t,x,\omega)\bigr)
-
\bigl(\nabla_x\varphi_k(t,x,\xi')+
\rho'\nabla_x\varphi_{j'}(t,x,\omega')\bigr)
\bigr|\ge c\,2^k\delta(l)\ge c\,2^{\frac{3k}4}\,.
$$
An easy integration by parts in $x$ using (1.10) and
the first part of (2.3) yields $(2.4)$.\qed
\enddemo

\proclaim{Lemma 2.2}
$$
\Bigl\|\,\int T^{0,\omega}(f_\mu,g_\omega)\,d\omega\,\Bigr\|_{L^2(dx\,dt)}
\le C\,
\|f_\mu\|_{L^2(\Bbb R^3)}\,\|g_\omega\|_{L^2(\Bbb R\times S^2)}\,.
$$
\endproclaim
\demo{Proof}
For fixed $j$ and fixed $(t,x)$,
the function $T^{0,\omega}(f_\mu,g_{j\omega})(t,x)$
vanishes unless $\omega$ is in a set of volume $\delta(0)^2=2^{-\frac k2}$.
Thus,
$$
\Bigl\|\, \int
T^{0,\omega}(f_\mu,g_\omega)\,d\omega\,\Bigr\|_{L^2(dx\,dt)}\le
\sum_j 2^{-\frac k4}
\,\left\|T^{0,\omega}(f_\mu,g_{j\omega})\right\|_{L^2(dx\,dt\,d\omega)}\,.
$$
Because of (2.3), the operator
$$
Af(x)=
\int e^{i\varphi_k(t,x,\xi)}
a_k(t,x,\xi)\,q^0_{kj}(t,x,\xi,\omega)
\,\widehat{f}(\xi)\,d\xi
\tag2.5$$
has $L^2\to L^2$ norm, for each fixed $t$,
less than $C\,\delta(0)=C\,2^{-\frac k4}$, with $C$ independent of $t$.
For the obstacle problem, where, for all $k$,  $g_k$ equal a fixed smooth
metric $g$, this just follows
from standard $L^2$ estimates for Fourier integral operators.  The general case where
there is a $k$-dependence also follows from standard $L^2$ estimates along with (1.6).
(See \cite{8}, \cite{9}.)

The aforementioned bounds for $Af$ immediately yield
$$
\left\|\,T^{0,\omega}(f_\mu,g_{j\omega})\right\|_{L^2(dx\,dt)}\le
C\,2^{-\frac k4}\,\|f_\mu\|_{L^2(\R^3)}\,
\|\widehat{g}_{j\omega}\|_{L^1(\R)}\le
C\,2^{-\frac k4}\,2^{\frac j2}\,
\|f_\mu\|_{L^2(\R^3)}\,\|g_{j\omega}\|_{L^2(\R)}\,.
$$
The lemma now follows since
$$
\sum_{j\le k}2^{(j-k)/2}\|g_{j\omega}\|_{L^2(\R\times S^2)}
\le\|g_\omega\|_{L^2(\R\times S^2)}\,.\qed
$$
\enddemo

\head
3. Null form Estimates
\endhead

\noindent
In this section, we show that, for each fixed $k$, $\mu$, and $l\ge 1$, the
following holds:
$$
\Bigl\|
\int T^{l,\omega}(f_\mu,g_\omega)\,d\omega\Bigr\|_{L^2(dx\,dt)}\le
C\,\delta(l)^{\frac 14}\,|\log\delta(l)|^{\frac 12}
\,\|f_\mu\|_{L^2(\Bbb R^3)}\,
\|g_\omega\|_{L^2(\Bbb R\times S^2)}\,,
\tag3.1$$
with constant $C$ independent of $k$, $l$, and $\mu$.
Together with Lemma 2.1 and Lemma 2.2, this implies estimate $(2.1)$,
after summing over $l$,
which in turn implies the desired estimate $(1.2)$.

We establish $(3.1)$ by splitting the operator $T^{l,\omega}$
into two pieces.
Let
$$
\aligned
T_1^{l,\omega}(f,g)&=\sum_{\{j\, :\, 2^j>2^{\frac k2}\delta(l)^{-1}\}}
T^{l,\omega}(f,g_j)\, ,
\\
\\
T_2^{l,\omega}(f,g)&=\sum_{\{j\, :\, 2^j\le 2^{\frac k2}\,\delta(l)^{-1}\}}
T^{l,\omega}(f,g_j)
\,,
\,.
\endaligned
\tag3.2$$

For the operator $T^{l,\omega}_1$,
note that $2^{-\frac j2}\le \delta(l)$ for the indices arising, since
$2^{\frac k2}\ge \delta(l)^{-1}$.
Hence, by the second part of (1.9) and the definition of $q_{kj}^l$,
the symbol of $T^{l,\omega}_1$ vanishes unless
$$
\text{angle}\,\bigl(\,d_x\varphi_k(t,x,\xi)\,,\,
d_x\varphi_k(t,x,\omega)\,\bigr)
\le C\,\delta(l)\,.
$$
Since the map
$\xi\rightarrow d_x\varphi_k(t,x,\xi)/|d_x\varphi_k(t,x,\xi)|$
is a $C^1$ diffeomorphism of the unit sphere, with uniform bounds over
$t,x,$ and $k$, for $t$ small,
it follows that the integrand vanishes unless
$\bigl|\xi^\mu-\omega\bigr|\le C\,\delta(l)\,.$
Consequently, by the Schwarz inequality
$$
\Bigl\|\,
\int T_1^{l,\omega}(f_\mu,g_\omega)\,d\omega\,\Bigr\|_{L^2(dx\,dt)}\le
\delta(l)\,
\bigl\|\,
T_1^{l,\omega}(f_\mu,g_\omega)\,\bigr\|_{L^2(dx\,dt\,d\omega)}\,.
$$
For the piece $T^{l,\omega}_1$, the estimate $(3.1)$ is thus
implied by the following
\proclaim{Theorem 3.1}
The following holds, with $C$ independent of $l\,,\,\omega\,,\,k\,,$
$$
\bigl\|\,T_1^{l,\omega}(f,g)\bigr\|_{L^2(dx\,dt)}
\le\,C\,\delta(l)^{-\frac 34}\,
|\log\delta(l)|^{\frac 12}\,
\,\|f\|_{L^2(\Bbb R^3)}\,\|g\|_{L^2(\Bbb R)}\,.
$$
\endproclaim

\smallskip

We postpone the proof of Theorem 3.1, and first establish the somewhat 
easier

\proclaim{Theorem 3.2}
The following holds, with $C$ independent of $l\,,\,\mu\,,\,k\,,$
$$
\Bigl\|
\int T_2^{l,\omega}(f_\mu,g_\omega)\,d\omega\Bigr\|_{L^2(dx\,dt)}\le
\,C\,\delta(l)\;|\log\delta(l)|
\;\,\|f_\mu\|_{L^2(\Bbb R^3)}\,\|g_\omega\|_{L^2(\Bbb R\times S^2)}\,.
$$
\endproclaim

\demo{Proof}
We split the sum over $j$ in (3.2) into three distinct cases: $2^j\le \delta(l)^{-2}$,
$\delta(l)^{-2}< 2^j\le 2^{k/2}$, and $2^{k/2}<2^j\le \delta(l)^{-2}2^{k/2}$.

\demo{\bf{Case 1: $2^j\le \delta(l)^{-2}$}}
In this case the index $j$ runs over $O(|\log\delta(l)|)$ values. Also,
for fixed $j$ and fixed $(t,x)$, the integrand vanishes unless
$\omega$ lies in a set in $S^2$ of area $\delta(l)^2$. Hence, by the Schwarz
inequality, it
suffices to establish the following estimate, uniformly in $j$
and $\omega$:
$$
\bigl\|
T^{l,\omega}(f,g_j)\bigr\|_{L^2(dx\,dt)}\le
C\,\|f\|_{L^2(\R^3)}\,\|g_j\|_{L^2(\R)}\,.
\tag3.3$$
We now write $\widehat f(\xi)=\sum_\nu \widehat{f}_\nu(\xi)$
where $\widehat{f}_\nu$ is supported in a cone
of angle $2^{-\frac k2}$ about a unit vector
$\xi^\nu$. The estimate $(3.3)$ is a result of the following
$$
\left|\,\int T^{l,\omega}(f_\nu,g_j)\,
\overline{T^{l,\omega}(f_{\nu'},g_j)}\,
dt\,dx\,\right|\le C\,\bigl(1+2^{\frac
k2}|\,\xi^\nu-\xi^{\nu'}|\,\bigr)^{-N}\,
\|f_\nu\|_{L^2(\Bbb R^3)}\,\|f_{\nu'}\|_{L^2(\Bbb R^3)}\,
\|g_j\|_{L^2(\Bbb R)}^2\,.
\tag3.4$$
The $(t,x)$ integrand in $(3.4)$ is dominated by
$$
\multline
\delta(l)^2\,
\left|\int
e^{i\varphi_k(t,x,\xi)-i\varphi_k(t,x,\xi')}a_{\nu,\nu'}(t,x,\xi,\xi')
\,\widehat{f}_\nu(\xi)\,\widehat{f}_{\nu'}(\xi')\,d\xi\,d\xi'\right|
\\
\times
\left|\int e^{i\rho\varphi_j(t,x,\omega)}\,a_j(t,x,\rho\omega)\,
\widehat{g}_j(\rho)\,d\rho\,\right|^2\,,
\endmultline
\tag3.5$$
where, by $(1.10)$ and $(2.3)$,
$$
\left|\partial_{t,x}^\beta\partial_{\xi,\xi'}^\alpha
\langle\xi^\nu,\partial_\xi\rangle^m\langle\xi^{\nu'},\partial_{\xi'}\rangle
^{m'}
a_{\nu,\nu'}(t,x,\xi,\xi')\right|\le
C\,2^{\frac k2 (|\beta|-|\alpha|)-k(m+m')}\,.
\tag3.6$$
Since $\rho\le\delta(l)^{-2}\le 2^{\frac k2}\,,$ the operator
$2^{-\frac k2}\partial_{x}$ applied to the expression inside the absolute
value sign in $(3.5)$ leads to an expression
of the same form. Furthermore, on the $(\xi,\xi')$ support of the
symbol in $(3.5)$, the following holds,
$$
2^{-\frac k2}
\bigl|\partial_x\varphi_k(t,x,\xi)-\partial_x\varphi_k(t,x,\xi')\bigr|\ge
c\,2^{\frac k2}\,|\,\xi^\nu-\xi^{\nu'}|\,.
$$
Integration by parts in $x$ now bounds the left hand side of $(3.4)$ by
$$
\multline
\frac{\delta(l)^2}{\bigl(1+2^{\frac k2}|\,\xi^\nu-\xi^{\nu'}|\,\bigr)^N}\,
\int\left|\,
\int e^{i\varphi_k(t,x,\xi)-i\varphi_k(t,x,\xi')}a_{\nu,\nu'}(t,x,\xi,\xi')
\,\widehat{f}_\nu(\xi)\,\widehat{f}_{\nu'}(\xi')\,d\xi\,d\xi'\,\right|\\
\times
\left|\,\int e^{i\rho\varphi_j(t,x,\omega)}\,a_j(t,x,\rho\omega)\,
\widehat{g}_j(\rho)\,d\rho\,\right|^2\,dt\,dx\,,
\endmultline
$$
where $a_{\nu,\nu'}(t,x,\xi,\xi')$ is a symbol satisfying the same
estimates $(3.6)$.
Next, following \cite{6}, we replace the phase $\varphi_k(t,x,\xi)$
by $\langle\nabla_\xi\varphi_k(t,x,\xi^\nu),\xi\rangle$, modulo an error
that is absorbed into the symbol, and similarly for $\varphi_k(t,x,\xi')\,.$
The left hand side of $(3.4)$ is thus bounded by
$$
\frac{\delta(l)^2}{\bigl(1+2^{\frac k2}|\,\xi^\nu-\xi^{\nu'}|\,\bigr)^N}\,
\int f_\nu^*\bigl(\nabla_\xi\varphi_k(t,x,\xi^\nu)\bigr)
\,f_{\nu'}^*\bigl(\nabla_\xi\varphi_k(t,x,\xi^{\nu'})\bigr)
\,g_j^*\bigl(\varphi_j(t,x,\omega)\bigr)^2\,dt\,dx\,,
\tag3.7$$
where
$$
f_\nu^*(y)=2^{2k}\int
\bigl(1+2^{\frac k2}|y-z|+2^k|\langle\xi^\nu,y-z\rangle|\,\bigr)^{-4}\,
|f_\nu(z)|\,dz\,,
$$
and
$$
g_j^*(s)=\int \bigl(1+2^j|s-r|\,\bigr)^{-2}\,|g_j(r)|\,ds\,,
$$
hence
$$
\|f^*_\nu\|_{L^2(dy)}\le C\|f_\nu\|_{L^2(\Bbb R^3)}\,,
\quad
\|g^*_j\|_{L^2(ds)}\le C\|g_j\|_{L^2(\Bbb R)}\,.
$$
The change of variables
$(t,x)\rightarrow\bigl(\varphi_j(t,x,\omega),\nabla_\xi\varphi_k(t,x,\xi^\nu)
\bigr)$
has Jacobian comparable to $\delta(l)^2$. An application of the Schwarz
inequality to $(3.7)$ thus yields $(3.4)$.
\enddemo

\demo{\bf{Case 2: $\delta(l)^{-2}< 2^j\le 2^{\frac k2}$}}
Consider the operator ${\tilde T}^{l,\omega}$ obtained by replacing
$q_{jk}^l(t,x,\xi,\omega)$ in equation $(2.2)$ by
$q_{jk}^l(t,x,\xi,\omega)-q_{kk}^l(t,x,\xi,\omega)\,.$
The proof of the previous case, together with the second set
of estimates in $(2.3)$, shows that for $2^j\le 2^{\frac k2}$
the following holds,
$$
\Bigl\|\,
\int {\tilde
T}^{l,\omega}(f_\mu,g_{j\omega})\,d\omega\,\Bigr\|_{L^2(dx\,dt)}
\le C\,2^{-\frac j2}\,
\,\|f_\mu\|_{L^2(\Bbb R^3)}\,\|g_{j\omega}\|_{L^2(\Bbb R\times S^2)}\,.
$$
Applying the Schwarz inequality over $j$ such that $2^j\ge\delta(l)^{-2}$
yields the estimate of Theorem 3.2 for this case
if $T^{l,\omega}$ is replaced by
${\tilde T}^{l,\omega}$. It thus remains to establish the same estimate
for the term
$$\multline
S^{l,\omega}(f_\mu,g_\omega)=
\\
\left(\;\int e^{i\varphi_k(t,x,\xi)}\,
a_k(t,x,\xi)\,q^l_{kk}(t,x,\xi,\omega)
\,\widehat{f}_\mu(\xi)\,d\xi\right)\,
\sum_j
\int
e^{i\rho\varphi_j(t,x,\omega)}\,a_j(t,x,\rho\omega)
\,\widehat{g}_j(\rho)\,d\rho\,.
\endmultline
$$
The $\xi$-integrand vanishes unless $|\omega-\xi^\mu|\le \delta(l)$,
hence
$$
\Bigl\|\,
\int S^{l,\omega}(f_\mu,g_{\omega})\,d\omega\,\Bigr\|_{L^2(dx\,dt)}\le
C\,
\delta(l)\,
\bigl\|
S^{l,\omega}(f_\mu,g_{\omega})\bigr\|_{L^2(dx\,dt\,d\omega)}\,.
$$
The proof of estimate $(3.4)$ establishes the following bound,
$$
\multline
\Bigl|\,\int S^{l,\omega}(f_\nu,g)\,
\overline{S^{l,\omega}(f_{\nu'},g)}\,
dt\,dx\,\Bigr|\\
\le C\,\bigl(1+2^{\frac k2}|\,\xi^\nu-\xi^{\nu'}|\,\bigr)^{-N}\,
\|f_\nu\|_{L^2(\Bbb R^3)}\,\|f_{\nu'}\|_{L^2(\Bbb R^3)}\,
\bigl(\,\sup_y\|Pg(\cdot,y)\|_{L^2(ds)}\,\bigr)\,
\bigl(\,\sup_{y'}\|P'g(\cdot,y')\|_{L^2(ds')}\bigr)\,,
\endmultline
$$
where $Pg$ is an operator of the form
$$
Pg=
\sum_j\int
e^{i\rho\varphi_j(t,x,\omega)}\,a_j(t,x,\rho\omega)
\,\widehat{g}_j(\rho)\,d\rho\,,
$$
written in the new coordinates
$$
(s,y)=\bigl(\varphi_k(t,x,\omega),\nabla_\xi\varphi_k(t,x,\xi^\nu)\bigr)\,,
$$
and $P'g$ is the same form with $\nu$ replaced by $\nu'$.
Using $(1.9)$ we may write $Pg$ in the form
$$
Pg(s,y)=\sum_j\int e^{is\rho}\,a_j(s,y,\rho)\,\widehat{g}_j(\rho)\,d\rho\,,
$$
where the new symbol satisfies
$$
\bigl|\partial_s^m\partial_\rho^n a_j(s,y,\rho)\bigr|\le
C\,\bigl(2^{\frac j2}\delta(l)^{-2}\bigr)^m\,2^{-jn}\,,\quad m\le 1\,.
$$
A simple integration by parts establishes the following bound,
$$
\multline
\left|\,\int e^{is(\rho-\rho')}\,a_j(s,y,\rho)\,\widehat{g}_j(\rho)\,
\overline{a_{j'}(s,y,\rho')\,\widehat{g}_{j'}(\rho')}\,d\rho\,d\rho'\,ds
\right|
\\
\le C\,\Bigl(1+\delta(l)^2\,2^{\max(j,j')/2}\Bigr)^{-1}\,
\|g_j\|_{L^2(\R)}\,\|g_{j'}\|_{L^2(\R)}\,.
\endmultline
\tag3.8
$$
Summing over $j\,,\,j'$ such that $2^j\ge\delta(l)^{-2}\,,\,
2^{j'}\ge\delta(l)^{-2}\,,$ yields the following,
$$
\sup_y\|Pg\|_{L^2(ds)}\le C\,|\log\delta(l)|\,\|g\|_{L^2(\Bbb R)}\,,
$$
which completes the proof for the second case.
\enddemo

\demo{\bf{Case 3: $2^{\frac k2}< 2^j \le 2^{\frac k2}\delta(l)^{-1}$}}
There are $O(|\log\delta(l)|)$ terms $j$, so as in the first case it
suffices
to establish the estimate $(3.4)$ uniformly over $j$. Let
$$
v(t,x)=\frac{1}{|\,\xi^\nu-\xi^{\nu'}|}\times\left(\,
\text{projection of}\,\;
\nabla_x\varphi_k(t,x,\xi^\nu)-\nabla_x\varphi_k(t,x,\xi^{\nu'})\,\;
\text{onto}\,\;\nabla_x\varphi_k(t,x,\xi^\nu)^\perp\,\right)\,.
$$
It follows from $(1.8)$ that
$$
\bigl|\,\partial_{t,x}^\alpha v(t,x)\bigr|\le C_\alpha\,2^{\frac k2
|\alpha|}\,.
$$
Next note that, if $q_{kj}^l(t,x,\xi,\omega)\,f_\nu(\xi)$ is nonzero,
then $|\xi^\nu-\omega|\le C\,\delta(l)\,.$ Also note that
$2^{-\frac j2}\le\delta(l)\,.$
The following is thus seen to hold by (1.8),
$$
\bigl|\,\partial_{t,x}^\alpha
\langle v(t,x),\nabla_x\varphi_j(t,x,\omega)\rangle\bigr|=
\bigl|\,\partial_{t,x}^\alpha
\langle v(t,x),\nabla_x\varphi_j(t,x,\omega)-\nabla_x\varphi_k(t,x,\xi^\nu)
\rangle\bigr|\le C_\alpha\,\delta(l)\,2^{\frac k2 |\alpha|}\,.
$$
Since $\rho\,\delta(l)\le 2^{\frac k2}\,,$ it follows that
for any $N$ one may write
$$
\bigl(\,2^{-\frac k2}\,\langle v(t,x),\nabla_x\rangle\,\bigr)^N
\left(\int e^{i\rho\varphi_j(t,x,\omega)}\,a_j(t,x,\rho\omega)\,
\widehat{g}_j(\rho)\,d\rho\,\right)
$$
as an expression of the same form as that in parentheses, 
but with a new symbol which satisfies the following estimates
$$
\bigl|
\partial_{t,x}^\alpha\partial_\rho^m \tilde a_j(t,x,\rho)\bigr|\le
C_{\alpha,m}\,2^{\frac k2 |\alpha|-mj}\,.
$$
These estimates imply the following bound,
$$
\left|\,\int e^{i\rho\varphi_j(t,x,\omega)}\,\tilde a_j(t,x,\rho\omega)\,
\widehat{g}_j(\rho)\,d\rho\,\right|\le C\,g^*_j(\varphi_j(t,x,\omega))\,.
$$
We next note that if $\widehat f_\nu(\xi)$ and 
$\widehat f_{\nu'}(\xi')$ are nonzero, and 
$|\xi^\nu-\xi^{\nu'}|\ge C\,2^{-\frac k2}$, then
$$
\langle v(t,x),\nabla_x\varphi_k(t,x,\xi)-\nabla_x\varphi_k(t,x,\xi')\rangle
\approx 2^k\,|\,\xi^\nu-\xi^{\nu'}|\,.
$$
The proof of estimate $(3.4)$ from the first case
now carries over to the
third case, where in establishing the estimate $(3.7)$ for the third
case, one integrates by parts using
$$\bigl(\langle
v(t,x),\nabla_x\varphi_k(t,x,\xi)-
\nabla_x\varphi_k(t,x,\xi')\rangle\bigr)^{-1}
\langle v(t,x),\nabla_x\rangle.$$

Since we have handled all three cases, the proof of Theorem 3.2 is complete. \qed

\enddemo
\enddemo

The proof of Theorem 3.1
rests on the following two lemmas estimating the gradients of the
phase function.

\proclaim{Lemma 3.3}
Let $\Phi(t,x,t',x',\xi)=\varphi(t,x,\xi)-\varphi(t',x',\xi)\,,$
where $\varphi=\varphi_k$ for some $k$.
Suppose that $\omega$ is a unit vector, and
let $\displaystyle\delta=\text{angle}(\omega,\xi)$.
If $|\xi|=2^k$, then for some $c>0$,
$$
\multline
2^{\frac k2}|\nabla_\xi\Phi(t,x,t',x',\xi)|
+
|\Phi(t,x,t',x',\xi)|
+
2^j|\Phi(t,x,t',x',\omega)|\\
\ge
c\,\Bigl(\,2^{\frac k2}\,|\nabla_\xi\Phi(t,x,t',x',\omega)|+
2^{\frac k2}\,\delta\,|\,t-t'|\,\Bigr)\,,
\endmultline
\tag3.9$$
for all $j$, $k$ and $\delta$ such that $j\le k$, and
$2^j\,\delta\ge 2^{\frac k2}\,.$
\endproclaim
\demo{Proof}
We have $2^k\ge 2^j\ge 2^{\frac k2}\delta^{-1}\,,$ so
introducing the new variables
$$
y_j(t,x)=\partial_{\xi_j}\varphi(t,x,\xi)\,, \qquad \mu=\xi/|\xi|,
$$
the left hand side of $(3.9)$ is larger than
$$
2^{\frac k2}|\,y-y'|
+
2^{\frac k2}\delta^{-1}\Bigl(\bigl|\langle\mu,y-y'\rangle\bigr|
+
|F(t,y)-F(t',y')|\Bigr)
\,,
\tag3.10
$$
where $F(t,y)$ is $\varphi(t,x,\omega)$ written in the coordinates $(t,y)$.
We begin by showing that the quantity $(3.10)$ is larger than
$$
2^{\frac k2}\bigr(\,|\,y-y'| + \delta\,|\,t-t'|\,\bigr)\,.
$$
To see this, note that the $C^1$ distance of $F(t,y)$ to
$\langle\mu,y\rangle$ is of size $\delta$, so that
$$
F(t',y)-F(t',y')=\langle\mu,y-y'\rangle+O\bigl(\,\delta|\,y-y'|\,\bigr)\,.
$$
Thus $(3.10)$ dominates
$$
2^{\frac k2}\bigr(\,|\,y-y'| +
\delta^{-1}|F(t,y)-F(t',y)|\,\bigr)\,.
$$
We will be done by establishing the following identity,
$$
\partial_t F(t,y)=\|d_x\varphi(t,x,\omega)\|_{\bold g}
-{\bold g}\bigr(d_x\varphi(t,x,\omega),d_x\varphi(t,x,\mu)\bigr)\,\big/\,
\|d_x\varphi(t,x,\mu)\|_{\bold g}\approx\delta^2\,,
$$
where $\bold g=\bold g_k\,.$
To see this, let $x=x(t,y)$ denote $x$ in the $(t,y)$ coordinates. Then
$\bigl(x,d_x\varphi(t,x,\mu)\bigr)$ is the backwards hamiltonian
curve through $(y,\mu)\,.$ Hamilton's equations thus yield
$$
\partial_t x_i=-\sum_{m=1}^3{\bold g}^{mi}(t,x)\partial_{x_m}
\varphi(t,x,\mu)\,\big/\,
\|d_x\varphi(t,x,\mu)\|_{\bold g}\,.
$$
Thus,
$$
\matrix\format\r\;&\c&\;\l\\
\partial_t F(t,y) & = & \partial_t\varphi(t,x(t,y),\omega)\\
\\
{}& = & \partial_t\varphi(t,x,\omega)+\sum_{i=1}^3
\partial_{x_i}\varphi(t,x,\omega)
\partial_t x_i\\
\\
{} & = &
\|d_x\varphi(t,x,\omega)\|_{\bold g}
-{\bold g}\bigr(d_x\varphi(t,x,\omega),d_x\varphi(t,x,\mu)\bigr)\,\big/\,
\|d_x\varphi(t,x,\mu)\|_{\bold g}\,.
\endmatrix
$$
To finish the proof of the lemma,
let $f_j(t,y)$ denote $\partial_{\xi_j}\varphi(t,x,\omega)$ in the
$(t,y)$ coordinates. Then the $C^1$ distance of $f_j$ to $y_j$
is comparable to $\delta$, so that $|\partial_t f_j(t,y)|\le \delta\,.$
Consequently,
$$
\bigl|\nabla_\xi\Phi(t,x,t',x',\omega)\bigr|\le\sum_{j=1}^3
|f_j(t,y)-f_j(t',y')|\le C\bigl(\,|\,y-y'|+\delta\,|\,t-t'|\,\bigr)\,.\qed
$$
\enddemo

\proclaim{Lemma 3.4}
Let $\Phi(t,x,\xi,\xi')=\varphi(t,x,\xi)-\varphi(t,x,\xi')\,,$
where $\varphi=\varphi_k$ for some $k$.
Suppose that
$$
\left|\;\omega-\frac{\xi}{|\xi|}\;\right|
\quad,\quad
\left|\;\omega-\frac{\xi'}{|\xi'|}\;\right|\,\in\,
\bigl[\,C^{-1}\delta\,,\,C\,\delta\,\bigr]\,,
$$
that $|\xi|\,,\,|\xi'|\in\bigl[2^{k-1},2^{k+1}\bigr]$,
and that $|\rho|\,,\,|\rho'|\in \bigl[0, 2^k\bigr]$. Then for some
$c>0$, independent of $k,\delta$,
$$
\Bigl|\,\nabla_{t,x}\Phi(t,x,\xi,\xi')+
\nabla_{t,x}\Phi(t,x,\rho\omega,\rho'\omega)\,\Bigr|
\ge
c\,\Bigl(\,
2^k\,\delta\,\times\text{\rm angle}(\xi,\xi')
+\delta^2\,|\rho-\rho'|\,
\Bigr)\,.
\tag3.11$$
\endproclaim
\demo{Proof}
Let $w=\nabla_x\varphi(t,x,\xi)\,,\;w'=\nabla_x\varphi(t,x,\xi')\,,
$ and $\mu=\nabla_x\varphi(t,x,\omega)\,.$
Also let $\alpha=\rho-\rho'\,.$
The conditions of the statement imply that the angle of
$w$ or $w'$ to $\mu$ is comparable to $\delta\,.$

By the eikonal equations,
the left hand side of $(3.11)$ dominates
$$
\bigl| w-w'+\alpha\mu\bigr|+\bigl|\,\|w\|-\|w'\|+\alpha\|\mu\|\,\bigr|
\tag3.12
$$
where $\|\cdot\|$ denotes the norm in the metric $\bold g_k(t,x)$.
We consider the case $\alpha\ge 0\,;$ the case $\alpha\le 0$ follows
by symmetry
upon exchanging $\xi$ and $\xi'\,.$ Also, by scaling $\alpha$, we may
assume that $\|\mu\|=1$. The quantity $(3.12)$ then dominates
$$
\matrix\format\r\;&\c&\;\l\\
\|w\|+\alpha-\|w+\alpha\mu\|\,& \ge & \,
c\,2^{-k}
\Bigl(\bigl(\,\|w\|+\alpha\bigr)^2-\|w+\alpha\mu\|^2\Bigr)\\
\\
{} & = & c\,2^{-k}\,\|w\|\,\alpha\,
\left\|\,\displaystyle\frac{w}{\|w\|}-\mu\,\right\|^2\\
\\
{} & \ge & c\,\delta^2\,\alpha\,.
\endmatrix
$$
We next observe that $(3.12)$ dominates the following quantity
(recall that $\|\mu\|=1$)
$$
2^k\,\left\|\,\frac{w}{\|w\|}-\mu
-r\left(\,\frac{w'}{\|w'\|}-\mu\,\right)\right\|\,,
\tag3.13$$
where $r=\|w'\|/\|w\|\in [c,c^{-1}]\,.$ By making a linear transformation,
we may replace the $\bold g$ norm $\|\cdot\|$ by the Euclidean
norm $|\cdot|$, and assume that
$$
\mu = (1,0,0)\,,\qquad
\frac{w}{|w|} = (\cos\theta,\sin\theta,0)\,,\qquad
\frac{w'}{|w'|} = (\sqrt{1-z^2}\,\cos\gamma,\sqrt{1-z^2}\,\sin\gamma,z)\,,
$$
where $\theta,\gamma,z$ are small. The quantity $(3.13)$ is then comparable
to
$$
2^k\,\Bigl(\,\bigl|\,1-\cos\theta-r\,(1-\sqrt{1-z^2}\,\cos\gamma)\,\bigr|+
\bigl|\,\sin\theta-r\sqrt{1-z^2}\,\sin\gamma\,\bigr|+|z|\,\Bigr)\,,
$$
which in turn is comparable to
$$
2^k\,\Bigl(\,\bigl|\,1-\cos\theta-r\,(1-\cos\gamma)\,\bigr|+
\bigl|\,\sin\theta-r\sin\gamma\,\bigr|+|z|\,\Bigr)\,.
$$
By the half angle formula, this equals
$$
2^k
\Bigl(\,\bigl|\,\sin\theta\,\tan{\theta\over 2}-2r\sin^2{\gamma\over 2}\,
\bigr|
+\bigl|\,\sin\theta-2r\sin{\gamma\over 2}\,\cos{\gamma\over 2}\,\bigr|+|z|
\,\Bigr)\,.
$$
Since $\cos(\gamma/2)\approx 1$, this in turn dominates
$$
2^k
\Bigl(\,|\sin\theta|\,\bigl|\,\tan{\theta\over 2}-\tan{\gamma\over 2}
\,\bigr|+|z|\,
\Bigr)
\ge c\,2^k\,\bigl(\,\delta\,|\,\theta-\gamma|+|z|\,\bigr)
\ge c\,2^k\,\delta\,\times
\text{angle}(\xi,\xi')
\,.\qed
$$
\enddemo

\demo{Proof of Theorem 3.1}
To complete the proof, we make a further decomposition
$$
T_1^{l,\omega}=\sum_{\nu,j,s,z}T^{l,\omega}_{\nu,j,s,z}\,,
$$
where
$$
T^{l,\omega}_{\nu,j,s,z}F(t,x)=
\int e^{i\varphi_k(t,x,\xi)+i\rho\varphi_k(t,x,\omega)}
a^{l,\omega}_{\nu,j,s,z}(t,x,\xi,\rho)\,\widehat{F}(\xi,\rho)
\,d\xi\,d\rho\,.
\tag3.14$$
The index $\nu$
corresponds to a set
of unit vectors $\xi^\nu$ evenly spaced
by $c\,2^{-\frac k2}\,\delta(l)^{-1}\,,$ for some small
$c$ to be determined independent of $k$ and $l$.
The index $j$ runs over the integers such that
$$
2^k \;\ge\; 2^j \;\ge\; 2^{\frac k2}\,\delta(l)^{-1}\,.
$$
The indices $z$ and $s$ run over lattices such that
$$
2^{\frac k2}z\in\Bbb Z^3\,,\qquad
2^{\frac k2}\,\delta(l)\,s\in\Bbb Z\,.
$$
The symbol $a^{l,\omega}_{\nu,j,s,z}(t,x,\xi,\rho)$ is supported in the set
where $|\xi|\in [2^{k-1},2^{k+1}]\,,\quad\rho\in [2^{j-1},2^{j+1}]\,,$
and where
$$
\left|\,\frac{\xi}{|\xi|}-\xi^\nu\,\right|\le c\,2^{-\frac
k2}\,\delta(l)^{-1}
\,,\qquad
\bigl|\,\nabla_\xi\varphi_k(t,x,\omega)-z\,\bigr|\le 2\cdot2^{-\frac k2}
\,,\qquad
|\,t-s\,|\le 2\cdot2^{-\frac k2}\,\delta(l)^{-1}\,.
$$
The symbol furthermore satisfies the estimates
$$
\bigl|
\partial_{t,x}^\beta\partial_\xi^\alpha\partial_\rho^m
\langle\xi,\partial_\xi\rangle^i
a^{l,\omega}_{\nu,j,s,z}(t,x,\xi,\rho)\bigr|\le
C\,\delta(l)\,
2^{-jm+\frac k2(|\beta|-|\alpha|)}\,.
\tag3.15$$

A few remarks are in order here. First, 
as a result of (1.8) and (1.9), the function
$$
e^{i\rho\varphi_j(t,x,\omega)-i\rho\varphi_k(t,x,\omega)}
$$
satisfies
the symbol estimates $(3.15)$,
which allowed us to replace the phase $\varphi_j(t,x,\omega)$
by $\varphi_k(t,x,\omega)$ in formula $(3.14)$.
Next, since $\delta(l)\ge 2^{-\frac k4}$, 
$$
2^{-\frac j2}\le 2^{-\frac k8}\,\delta(l)\ll \delta(l)\, .
$$
It follows from $(1.9)$ and the definition of $q^l_{kj}$ that
the angle of $\nabla_x\varphi_k(t,x,\omega)$ to $\nabla_x\varphi_k(t,x,\xi)$
is comparable to $\delta(l)$, hence that the angle of $\omega$ to $\xi$
is comparable to $\delta(l)$. By making the number $c$ above small,
it follows that the angle of $\omega$ to $\xi^\nu$ is comparable to
$\delta(l)$.

We begin by showing that
$$
\sup_{s',z'}\;\sum_{s,z}\;
\bigl\|\,
T^{l,\omega}_{\nu,j,s,z}\bigl(\,T^{l,\omega}_{\nu,j,s',z'}\,\bigr)^*
\,\bigr\|^{\frac 12}
\le C\,\delta(l)^{-1}\,.
\tag3.16$$
This will follow from showing that
$$
\bigl\|\,
T^{l,\omega}_{\nu,j,s,z}\bigl(\,T^{l,\omega}_{\nu,j,s',z'}\,\bigr)^*
\,\bigr\|
\;\le\;
\frac{C_N\,\delta(l)^{-2}}
{\bigl(\,1+2^{\frac k2}|\,z-z'|+2^{\frac
k2}\,\delta(l)\,|\,s-s'|\,\bigr)^N}\,.
\tag3.17$$
To establish $(3.17)$, we express
$T^{l,\omega}_{\nu,j,s,z}\bigl(\,T^{l,\omega}_{\nu,j,s',z'}\,\bigr)^*$
as an integral kernel of the form
$$
\multline
K(t,x;t',x')=
\\
\int
e^{i\varphi_k(t,x,\xi)-i\varphi_k(t',x',\xi)
+i\rho(\varphi_k(t,x,\omega)-\varphi_k(t',x',\omega))}
a^{l,\omega}_{\nu,j,s,z}(t,x,\xi,\rho)
\,
\overline{a}^{l,\omega}_{\nu,j,s',z'}(t',x',\xi,\rho)\,
d\xi\,d\rho\,.
\endmultline
$$
Integration by parts in $\xi$ and $\rho$, together with
the estimates $(3.15)$ and the support conditions, shows that
the kernel $|K(t,x;t',x')|$ is bounded by
$$
\int_{R_{\nu,j}}
\frac{C_N\,\delta(l)^2}
{\bigl(\,1+
2^{\frac k2}\bigl|\,\nabla_\xi\Phi(t,x,t',x',\xi)\,\bigr|
+
\bigl|\,\Phi(t,x,t',x',\xi)\,\bigr|
+
2^j\bigl|\,\Phi(t,x,t',x',\omega)\,\bigr|\,\bigr)^N}\,d\xi\,d\rho
\,,
$$
where $\Phi$ is as in Lemma 3.3.
For each $\xi$ in the domain of integration, the change of variables
$(t,x)\rightarrow \bigl(\nabla_\xi\varphi_k(t,x,\xi),\varphi_k(t,x,\omega)
\bigr)$
has Jacobian comparable to $\delta(l)^2$;
consequently, by Schur's Lemma and Lemma 3.3, for each fixed $\xi$ and
$\rho$
the integrand is an operator on $L^2(dt\,dx)$ with norm bounded by
$$
\frac{C_N\,2^{-2k-j}}
{\bigl(\,1+2^{\frac k2}|\,z-z'|+2^{\frac
k2}\,\delta(l)\,|\,s-s'|\,\bigr)^N}\,.
$$
The volume of $R_{\nu,j}$ is comparable to $2^{2k+j}\,\delta(l)^{-2}\,,$
and the estimate $(3.17)$ follows.

We next establish the following estimate
$$
\sup_{j',\nu'}\;\sum_{j,\nu}\;
\bigl\|
\,\bigl(T^{l,\omega}_{\nu,j,s,z}\,\bigr)^*\,
T^{l,\omega}_{\nu',j',s,z}\,\bigr\|^{\frac 12}
\;\le\;
C\,\delta(l)^{-\frac 12}\,|\log\delta(l)|\,.
\tag3.18$$
This will follow from showing that
$$
\bigl\|
\,\bigl(T^{l,\omega}_{\nu,j,s,z}\,\bigr)^*\,
T^{l,\omega}_{\nu',j',s,z}\,\bigr\|
\;\le\;
\frac{C_N\,\delta(l)^{-1}}
{\bigl(\,1+2^{-\frac k2}\,\delta(l)^2\,\bigl|\,2^j-2^{j'}\bigr|
+2^{\frac k2}\,\delta(l)\,|\,\xi^\nu-\xi^{\nu'}|\,
\bigr)^N}\,.
\tag3.19$$
provided that $|\,j-j'|\ge 3\,.$ For $|\,j-j'|\le 2\,,$ the estimate
holds as if $j=j'\,.$
That $(3.18)$ is a result of $(3.19)$ follows from that fact that
$$
\sum_j\;
\frac{1}
{1+2^{-\frac k2}\,\delta(l)^2\,|\,2^j-2^{j'}|}\le C\,|\log\delta(l)|\,,
$$
where the sum is over $j$ such that
$2^j\ge 2^{\frac k2}\,\delta(l)^{-1}$.

To establish $(3.19)$, we note that
$\bigl(\,T^{l,\omega}_{\nu,j,s,z}\,\bigr)^*T^{l,\omega}_{\nu',j',s,z}$
has an integral kernel of the form
$$
\multline
K(\xi,\rho;\xi',\rho')=
\\
\int
e^{i\varphi_k(t,x,\xi)-i\varphi_k(t,x,\xi')
+i\varphi_k(t,x,\rho\omega)-i\varphi_k(t,x,\rho'\omega)}
\overline{a}^{l,\omega}_{\nu,j,s,z}(t,x,\xi,\rho)
\,
a^{l,\omega}_{\nu',j',s,z}(t,x,\xi',\rho')\,
dt\,dx\,.
\endmultline
$$
Integration by parts in $(t,x)$ yields the following bound,
$$
|K(\xi,\rho;\xi',\rho')|\;\le\;
\int_{R_{s,z}}
\frac{C_N\,\delta(l)^2}
{\Bigl(\,1+2^{-\frac k2}\Bigl|\,\nabla_{t,x}\Phi(t,x,\xi,\xi')+
\nabla_{t,x}\Phi(t,x,\rho\omega,\rho'\omega)\,\Bigr|\,\Bigr)^N}\,dt\,dx\,,
$$
where $\Phi$ is as in Lemma 3.4, and $R_{s,z}$ is a set of
volume $2^{-2k}\delta(l)^{-1}\,.$ The change of variables
$(\xi,\rho)\rightarrow \nabla_{t,x}\bigl(\varphi_k(t,x,\xi)+
\varphi_k(t,x,\rho\omega)\bigr)$ has, for each fixed $(t,x)$, Jacobian
factor
comparable to $\delta(l)^2\,.$
The estimate $(3.19)$ now follows from Schur's Lemma and Lemma 3.4.

To conclude the proof of Theorem 3.1,
we split $T^{l,\omega}$ into a finite number of pieces so that
we may assume that
$$
\bigl(\,T^{l,\omega}_{\nu,j,s,z}\,\bigr)^*T^{l,\omega}_{\nu',j',s',z'}=0
$$
unless $z=z'$ and $s=s$, and
$$
T^{l,\omega}_{\nu',j',s',z'}\bigl(\,T^{l,\omega}_{\nu,j,s,z}\bigr)^*=0
$$
unless $\nu=\nu'$ and $j=j'$.
We now consider an arbitrary finite truncation of
the following sum to $M$ elements
$$
T^{l,\omega}=\sum_{\nu,j,s,z}T^{l,\omega}_{\nu,j,s,z}
$$
The proof of the Cotlar-Stein Lemma yields the following,
$$
\multline
\bigl\|\,T^{l,\omega}\,\bigr\|^{2N}\le
C\,\delta(l)^{-1}\,\sum\;
\bigl\|\,\bigl(\,T^{l,\omega}_{\nu_1,j_1,s_1,z_1}\,\bigr)^*
T^{l,\omega}_{\nu_2,j_2,s_1,z_1}\,\bigr\|^{\frac 12}\,
\bigl\|\,T^{l,\omega}_{\nu_2,j_2,s_1z_1}
\bigl(\,T^{l,\omega}_{\nu_2,j_2,s_2z_2}\,\bigr)^*\bigr\|^{\frac 12}\,\\
\bigl\|\,\bigl(\,T^{l,\omega}_{\nu_2,j_2,s_2z_2}\,\bigr)^*
T^{l,\omega}_{\nu_3,j_3,s_2z_2}\,\bigr\|^{\frac 12}
\cdots\,
\bigl\|\,\bigl(\,T^{l,\omega}_{\nu_N,j_N,s_N,z_N}\bigr)^*
T^{l,\omega}_{\nu_{N+1},j_{N+1},s_N,z_N}\bigr\|^{\frac 12}
\endmultline
$$
and by estimates $(3.16)$ and $(3.18)$ this implies
$$
\bigl\|\,T^{l,\omega}\,\bigr\|^{2N}\le
M\,C^{2N}\,\delta(l)^{-\frac{3N}2}\,|\log\delta(l)|^N\,.
$$
Letting $N\rightarrow\infty$ completes the proof of Theorem 3.1.\qed

\enddemo

\head
4. Null form Estimates for the Wave Equation on Geodesically Concave
Manifolds
\endhead

In this section we work locally on a three-dimensional
Riemannian manifold $\Omega$ with metric $\bold g$ and
with smooth boundary $\partial\Omega$, such
that $\Omega$ is strictly geodesically concave with respect to $\bold g$.
The typical example is $\Omega$ the complement in $\R^3$ of a strictly
convex open set, with the Euclidean metric understood.
By the Cauchy problem on $\Omega$ with Dirichlet condition we
understand the following system
$$
\cases
\partial_t^2 u(t,x)=\Delta_{\bold g} u(t,x)+F(t,x)\,,
\\
u(t,x) = 0\quad\text{if}\quad x\in\partial\Omega\,,
\\
u(0,x)=u_0(x), \quad \partial_t u(0,x)=u_1(x).
\endcases
$$
We work in a local coordinate patch centered at the origin such that
$\Omega$ is defined by $x_3\ge 0$. For $k=1,2$ we set
$$
H^k_D(\Omega)=\bigl\{f\in H^k(\Omega)\,:\,f|_{\partial\Omega}=0\,\bigr\}
\,,
$$
where $H^k(\Omega)$ is the space of restrictions of elements of $H^k(\R^3)$.
\proclaim{Theorem 4.1}
Suppose that $u$ and $v$ satisfy the Cauchy problem on $\Omega$
with Dirichlet condition, with respective data
$$
u_0,v_0\in H^2_D(\Omega)\,,\quad
u_1,v_1\in H^1_D(\Omega)\,,\quad
F,G,DF,DG\in L^1_t(\,[-\delta,\delta];L^2(\Omega))\,.
$$
Suppose also that the data vanishes for $|x|\ge \delta$, where $\delta>0$
is a constant depending on $\Omega$. Then the following hold,
for any of the null forms $Q$,
$$
\multline
\bigl\|DQ(du,dv)\bigr\|_{L^2_{t,x}(\,[-\delta,\delta]\times\Omega)}\le
C\,
\Bigl(\,\|u_0\|_{H^2_D(\Omega)}+\|u_1\|_{H^1_D(\Omega)}
+\sum_{|\alpha|\le 1}
\|D^\alpha F\|_{L^1_tL^2_x(\,[-\delta,\delta]\times\Omega)}\,\Bigr)\\
\times
\Bigl(\,\|v_0\|_{H^2_D(\Omega)}+\|v_1\|_{H^1_D(\Omega)}
+\sum_{|\alpha|\le 1}
\|D^\alpha G\|_{L^1_tL^2_x(\,[-\delta,\delta]\times\Omega)}\,\Bigr)
\,.
\endmultline
$$
$$
\multline
\bigl\|Q(du,dv)\bigr\|_{L^2_{t,x}(\,[-\delta,\delta]\times\Omega)}\le
C\,
\Bigl(\,\|u_0\|_{H^1_D(\Omega)}+\|u_1\|_{L^2(\Omega)}
+\|F\|_{L^1_tL^2_x(\,[-\delta,\delta]\times\Omega)}\,\Bigr)\\
\times
\Bigl(\,\|v_0\|_{H^2_D(\Omega)}+\|v_1\|_{H^1_D(\Omega)}
+\sum_{|\alpha|\le 1}
\|D^\alpha G\|_{L^1_tL^2_x(\,[-\delta,\delta]\times\Omega)}\,\Bigr)
\,.
\endmultline
$$
\endproclaim

Before proving this result, we should point out that it immediately
yields Theorem 1.1.  This just follows from the standard existence
argument given in \cite{2}.

\demo{Proof of Theorem 4.1}
For convenience, in this proof we refer to
the discussion in \cite{10} regarding
the parametrix for the Dirichlet problem; however, all of the results
used are due to Melrose and Taylor \cite{3}, \cite{4}, \cite {5}, and
Zworski \cite{13}.
Since we are working locally,
we may assume that $\Omega$ is a compact manifold, hence that the
Dirichlet Laplacian $-\Delta$ is strictly positive on $L^2\,.$

In the estimate for $DQ$, The terms where the $D$ acts on
the coefficients of $Q$ may be handled by energy estimates. Hence, by
symmetry we may replace 
$\|DQ(du,dv)\|_{L^2(\,[-\delta,\delta]\times\Omega)}$ by
$\|Q(d\partial u,dv)\|_{L^2(\,[-\delta,\delta]\times\Omega)}\,,$
where $\partial u$ is any space or time derivative of $u$.
The next step is to reduce Theorem 4.1 to
the following pair of estimates for the homogeneous problem,
$$
\matrix\format\r&\,\c\,&\l\\
\bigl\|Q(d\partial_x u,dv)\bigr\|_{L^2(\,[-\delta,\delta]\times\Omega)}&\le&
C\,
\Bigl(\,\|u_0\|_{H^2_D(\Omega)}+\|u_1\|_{H^1_D(\Omega)}
\,\Bigr)
\Bigl(\,\|v_0\|_{H^2_D(\Omega)}+\|v_1\|_{H^1_D(\Omega)}
\,\Bigr)\,,\\
\\
\bigl\|Q(du,dv)\bigr\|_{L^2(\,[-\delta,\delta]\times\Omega)}&\le&
C\,
\Bigl(\,\|u_0\|_{H^1_D(\Omega)}+\|u_1\|_{L^2(\Omega)}
\,\Bigr)
\Bigl(\,\|v_0\|_{H^2_D(\Omega)}+\|v_1\|_{H^1_D(\Omega)}
\,\Bigr)\,.
\endmatrix
\tag4.1
$$
To do this, we first reduce Theorem 4.1
to the case $G=0\,.$ To this end, we integrate by parts
to write the contribution to $v$ from $G$ as
$$
\multline
\int_0^t\frac{\sin\bigl((t-s)\sqrt{-\Delta}\,\bigr)}{\sqrt{-\Delta}}
\;G(s,x)\,ds\\
\matrix\format\c&\;\;\l\\
=&\displaystyle\cos\bigl(t\sqrt{-\Delta}\,\bigr)\Delta^{-1}\,G(0,x)
-\Delta^{-1}\,G(t,x)
+
\int_0^t\cos\bigl((t-s)\sqrt{-\Delta}\,\bigr)\Delta^{-1}
\,\partial_s G(s,x)\,ds
\\
\\
=&I+II+III\,,
\endmatrix
\endmultline
\tag4.2
$$
where $\Delta^{-1}$ denotes the inverse Laplacian on $\Omega$
with Dirichlet conditions, which maps $H^k(\Omega)$ to $H^{k+2}(\Omega)$
by elliptic regularity.

To handle $I$, we note that
$$
\|\Delta^{-1}G(0,\cdot\,)\|_{H^2_D(\Omega)}\;\le\;
C\,\|G(0,\cdot\,)\|_{L^2(\Omega)}\;\le\;
C\,\sum_{j\le 1}\|\partial_t^jG\|_{L^1_tL^2_x([-\delta,\delta]\times\Omega)}
\,.
$$
This term can thus be absorbed into the initial data $v_0\,.$

Next, let $\tilde v(t,x,s)=\cos\bigl((t-s)\sqrt{-\Delta}\,\bigr)\Delta^{-1}
\,\partial_s G(s,x)\,.$ 
Then $\tilde v(t,x,s)$ is a solution of the homogeneous wave equation
in $(t,x)$ for each $s$, with initial data satisfying
$$
\|\tilde v(0,\,\cdot\,,s)\|_{H^2_D(\Omega)}+
\|\partial_t\tilde v(0,\,\cdot\,,s)
\|_{H^1_D(\Omega)}\le
C\,\|\partial_s G(s,\cdot\,)\|_{L^2(\Omega)}\,.
$$
Note that the $t$-derivative of $II$ cancels the term in the
$t$-derivative of $III$ coming from the upper limit of integration.
Hence, we may write
$$
d(II+III)=
\int_0^t d\tilde v(t,x,s)\,ds\;+\;d_x(II)\,.
$$
Assuming that the second estimate of Theorem 4.1 holds in the case $G=0$,
we may bound
$$
\matrix\format\r&\;\c\;&\l\\
\displaystyle
\|Q(du,\int_0^t d\tilde v(\,\cdot,s)\,ds)\|_{L^2([-\delta,\delta]\times\Omega)}
&\le &
\displaystyle
\int_{-\delta}^\delta
\|Q(du,d\tilde v(\,\cdot,s))\|_{L^2([-\delta,\delta]\times\Omega)}\,ds\\
\\
&\le&
\displaystyle C
\Bigl(\,\|u_0\|_{H^1_D(\Omega)}+\|u_1\|_{L^2(\Omega)}+
\|F\|_{L^1_tL^2_x}
\,\Bigr)\,\|\partial_t G\|_{L^1_tL^2_x}\,.
\endmatrix
$$
The first estimate of the theorem is handled identically.

It remains to handle the term $d_x(II)\,.$ We do this by showing that
$$
\|d_x\,\Delta^{-1}G\|_{L^2_tL^\infty_x([-\delta,\delta]\times\Omega)}
\le C\sum_{|\alpha|\le 1}
\|D^\alpha G\|_{L^1_tL^2_x([-\delta,\delta]\times\Omega)}\,.
\tag4.3$$
Energy estimates show that $\|d\partial_x u\|_{L^\infty_tL^2_x}$
and $\|du\|_{L^\infty_tL^2_x}$ are bounded by the
appropriate norms of $u_0,u_1,$ and $F$, yielding the
desired estimate.

The proof of $(4.3)$ is based on the following estimate, which holds
globally on $\R^3$ for functions $f$ such that 
$\widehat f(\xi)\in L^1_{loc}\,,$
$$
\|f\|^2_{L^\infty(\R^3)}\le 
C\,\bigl\||D|f\bigr\|_{L^2(\R^3)}\,\bigl\|\Delta f\bigr\|_{L^2(\R^3)}\,.
$$ 
This estimate is verified by noting that it is dilation invariant,
so that one may reduce to the case 
$\bigl\|\Delta f\bigr\|_{L^2(\R^3)}=\bigl\||D|f\bigr\|_{L^2(\R^3)}=1\,,$ 
for which it follows easily by separately considering the low and
high frequencies of $f\,.$
We then bound
$$
\matrix\format\r&\;\c\;&\l\\
\displaystyle
\|d_x\,\Delta^{-1}G\|_{L^2_tL^\infty_x([-\delta,\delta]\times\Omega)}
&\le&
\displaystyle
C\,\int_{-\delta}^\delta \|G(t,\cdot)\|_{L^2(\Omega)}\,
\|G(t,\cdot)\|_{H^1(\Omega)}\,dt\\
\\
&\le&
\displaystyle
C\;\|G\|_{L^\infty_tL^2_x([-\delta,\delta]\times\Omega)}\,
\|G\|_{L^1_tH^1_x([-\delta,\delta]\times\Omega)}\\
\\
&\le& 
C\left(\sum_{|\alpha|\le 1}
\|D^\alpha G\|_{L^1_tL^2_x([-\delta,\delta]\times\Omega)}\right)^2\,,
\endmatrix
$$
which concludes the proof of $(4.3)$, and the reduction of the theorem
to the case $G=0\,.$

It remains to reduce Theorem 4.1 to the case $F=0\,.$ 
Consider the second estimate of the theorem. We note that
$$
d\int_0^t\frac{\sin\bigl((t-s)\sqrt{-\Delta}\,\bigr)}{\sqrt{-\Delta}}
\;F(s,x)\,ds=
\int_0^t\,d\left(\frac{\sin\bigl((t-s)\sqrt{-\Delta}\,\bigr)}{\sqrt{-\Delta}}
\;F(s,x)\right)\,ds\,,
$$ 
which reduces the second estimate to the case $F=0$; that is, the second
estimate of $(4.1)$.

As we have remarked previously, 
the first estimate of the theorem is reduced to considering 
$\|Q(d\partial u,dv)\|_{L^2(\,[-\delta,\delta]\times\Omega)}\,.$
To handle $Q(d\,\partial_t u,dv)$,
we note that $\partial_t u$ solves the Cauchy problem with data
in $H^1_D(\Omega)\times L^2(\Omega)$, with inhomogeneity in $L^1_tL^2_x$,
thus controlling $\|Q(d\,\partial_t u,dv)\|_{L^2}$
is reduced to the second estimate of Theorem 4.1, which we have
already reduced to $(4.1)$.

Next consider $Q(d\partial_x u,dv)$. 
We apply the identity $(4.2)$ with $G$ replaced by $F$, and
as before reduce to considering the term $Q(\partial_x^2(II),dv)\,.$
To bound the $L^2_{t,x}$ norm of this term, we note that
$$
\matrix\format\r&\;\c\;&\l\\
\bigl\|\,Q(\partial_x^2\Delta^{-1}F,dv)\,
\bigr\|^2_{L^2_{t,x}([-\delta,\delta]\times\Omega)}&\le&
\bigl\|\partial_x^2\Delta^{-1}F
\bigr\|^2_{L^2_tL^3_x([-\delta,\delta]\times\Omega)}
\,\|dv\|^2_{L^\infty_tL^6_x([-\delta,\delta]\times\Omega)}\\
\\
&\le&
C\,\|F\|_{L^\infty_tL^2_x([-\delta,\delta]\times\Omega)}\,
\|F\|_{L^1_tL^6_x([-\delta,\delta]\times\Omega)}\,
\|dv\|^2_{L^\infty_tL^6_x([-\delta,\delta]\times\Omega)}\\
\\
&\le&
C\,\left(\sum_{|\alpha|\le 1}
\|D^\alpha F\|_{L^1_tL^2_x([-\delta,\delta]\times\Omega)}\right)^2\,
\bigl(\,\|v_0\|_{H^1_D(\Omega)}+\|v_1\|_{L^2(\Omega)}\,\bigr)^2\,.
\endmatrix
$$
This concludes the reduction of Theorem 4.1 to the pair of estimates
$(4.1)$.

To establish the estimates $(4.1)$,
we note that, as discussed in \cite{10} immediately preceeding formulas
$(2.12)$ and $(2.24)$ of that paper,
for some $\delta$ as in the statement of the
theorem, the solution $v$ may be written, modulo smoothing
operators acting on the data, as a finite sum of terms of the form
$$
Tg(t,x)=\int e^{i\varphi^\pm(t,x,\xi)}\,a(t,x,\xi)\,
\widehat{g}(\xi)\,d\xi\,,
$$
where the phases are the solutions to the eikonal equation for some
smooth extension of the metric $\bold g$ to an open neighborhood of
the origin in $\R^3$, and the data $g\in H^2(\R^3)$ satisfies
$$
\|g\|_{H^2(\R^3)}\le
C\,\bigl(\,\|v_0\|_{H^2_D(\Omega)}+\|v_1\|_{H^1_D(\Omega)}\bigr)\,.
$$
The solution $u$ may be similarly written, with data $f$ belonging
respectively to $H^2(\R^3)$ or $H^1(\R^3)$, in the cases of the
two estimates $(4.1)$.
The amplitude $a(t,x,\xi)$, which is smooth in all variables and
vanishes for $|x|\ge C\,\delta\,,$ is of one of two types. 
Either it satisfies the modified
$S^0_{\frac 23,\frac 13}$ estimates $(1.4)$ of this paper, or it
satisfies the following estimates:
$$
\bigl|\,x_3^j\partial_{x_3}^k\,\langle\xi,\partial_\xi\rangle^N
\partial_{t,x_1,x_2}^\beta\partial_\xi^\alpha a(t,x,\xi)\,\bigr|
\le C_{j,k,N,\alpha,\beta}\,
\bigl(1+|\xi|\bigr)^{\frac 23(k-j-|\alpha|)+\frac 13|\beta|}\,.
\tag4.4$$
(We remark that in \cite{10} these estimates on the symbol were
shown to hold for $N=0\,;$ that the estimates hold
for general $N$ follows from the fact that these modified
estimates are preserved
under the equivalence of phase theorem of H\"ormander as seen,
for example, by the asymptotic formula for the transformed symbol,
and the fact that the symbol in our case is obtained by a change of phase
from the product of a standard symbol with cutoff functions that
satisfy (4.4).)

In either case, the operator $\partial_x T$ is an operator of the same
type, with a symbol of one higher order, hence the estimates $(4.1)$,
and consequently Theorem 4.1, are reduced to verifying the following
estimate
$$
\bigl\|Q(d\,T\!f,d\,T\!g)\bigr\|_{L^2_{t,x}([-\delta,\delta]\times\Omega)}\le
C\,\|f\|_{H^1(\Bbb R^3)}\|g\|_{H^2(\Bbb R^3)}\,,
\tag4.5$$
for $T$ an operator as above with a symbol satisfying either $(1.4)$
or $(4.4)$.

We remark that in \cite{10}, the Strichartz estimates were
shown to hold for both symbol types:
$$
\|Tf\|_{L^4_tL^4_x([-\delta,\delta]\times\Omega)}
\le C\,\|f\|_{H^{\frac 12}(\R^3)}\,.
$$

We first verify that the reductions of the second section of this paper
hold for symbols satisfying the estimates $(4.4)$. There are two places
where the arguments need to be modified.
The first is to verify that the estimate $(4.5)$
holds if, in the formula for $d\,T$, the $d$ acts on the symbol
$a(t,x,\xi)$.
Consider the term $d\,T\!f$, where the $d$
hits the symbol satisfying $(4.4)$. In this case, one
obtains an operator $Sf$ of the same form but with symbol of order
$\frac 23$. The resulting contribution to the left hand side of
$(4.5)$ is controlled by noting that
$$
\bigl\|\bigl(Sf)\,(d\,T\!g)\bigr\|_{L^2_{t,x}}\le
\|Sf\|_{L^6_tL^3_x}\,\|d\,T\!g\|_{L^\infty_tL^6_x}\le
C\|f\|_{H^1}\,\|g\|_{H^2}\,,
$$
where the last estimate for $Sf$ follows by interpolating the
following estimates
$$
\matrix\format\r\;&\c&\;\l\\
\|Sf\|_{L^4_tL^4_x} & \le & C\,\|f\|_{H^{\frac 76}(\R^3)}\,,\\
\\
\|Sf\|_{L^\infty_tL^2_x} & \le & C\,\|f\|_{H^{\frac 23}(\R^3)}\,.
\endmatrix
$$
Similarly one may bound
$$
\bigl\|(d\,T\!f)\,\bigl(Sg)\bigr\|_{L^2_{t,x}}\le
\|d\,T\!f\|_{L^\infty_tL^2_x}\,\|\,|D_x|^{\frac 56}Sg\|_{L^4_tL^4_x}\le
C\,\|f\|_{H^1}\,\|g\|_{H^2}\,.
$$

The other modification is to verify that the operator $(2.5)$
has norm of order $2^{-\frac k4}$, if now the symbol $a_k(t,x,\xi)$
satisfies $(4.4)$. This follows by expressing
$$
Af(x)=\int_0^{x_n}
2^{\frac 23 k}\,\bigl(\,1+2^{\frac 43 k}r^2\,\bigr)^{-1}\,
A_rf(x)\,dr
$$
where $A_r$ is the operator obtained by replacing $a_k(t,x,\xi)$ by
the symbol
$$
a_{k,r}(t,\overline x,\xi)=
2^{-\frac 23 k}\,\bigl(\,1+2^{\frac 43 k}r^2\,\bigr)\,
\partial_{x_3}a_k(t,\overline x,r,\xi)\,,\qquad
\overline x = (x_1,x_2)
$$
which satisfies, for each $r$, the estimates $(1.4)$, with
constants independent of $r$. One then has the bound
$$
\|Af\|_{L^2(\R^3_x)}\le \sup_r\,\|A_rf\|_{L^2(\R^3_x)}\le
C\,2^{-\frac k4}\|f\|_{L^2}\,,
$$
with, as before, the $2^{-\frac k4}=\delta(0)$ factor coming from (2.3).
This procedure of ``freezing the $x_3$ coefficient'' will be
used in subsequent steps.

We are thus reduced to establishing estimate $(3.1)$.
The above technique of freezing the $x_3$ coefficient reduces
to the case that the symbol $a_k(t,x,\xi)$ in formula $(2.2)$
satisfies the good estimates $(1.4)$, and the symbol $a_j(t,x,\rho\omega)$
satisfies the estimates $(4.4)$ above. (Note that one cannot freeze
the $x_3$ coefficient of $a_j(t,x,\rho\omega)$, since $\widehat g(\rho)$
is not localised to a dyadic interval.)

We next note that the proofs of Theorems 3.1 and 3.2
go through if $\widehat g(\rho)$ is supported in the
region where $\rho\le 2^{\frac {3k}4}$. This follows since, in this
case, we have $2^{\frac{2j}3}\le 2^{\frac k2}$, hence $\partial_x$
loses at most $2^{\frac k2}$ against the symbol $a_j(t,x,\rho\omega)$.
The only step in the proof that needs to be modified is to replace
the right hand side of $(3.8)$ by
$$
C\,\Bigl(1+\delta(l)^2\,2^{\max(j,j')/3}\Bigr)^{-1}\,
\|g_j\|_{L^2(\R)}\,\|g_{j'}\|_{L^2(\R)}\,,
$$
to reflect the $(\frac 23,\frac 23)$
estimates on $a_j(t,x,\rho\omega)$.

We thus assume that $\widehat g(\rho)$ is supported in the
region where $\rho\ge 2^{\frac {3k}4}$. Notice that
$\rho\ge 2^{\frac k2}\,\delta(l)^{-1}$, since
$\delta(l)\ge 2^{-\frac k4}$. Consequently
$T^{l,\omega}(f,g)=T^{l,\omega}_1(f,g)$. We will show that
$$
\bigl\|\,T^{l,\omega}(f,g)\bigr\|_{L^2(dx\,dt)}
\le\,C\,\delta(l)^{-\frac 34}\,
|\log\delta(l)|^{\frac 32}\,
\,\|f\|_{L^2(\Bbb R^3)}\,\|g\|_{L^2(\Bbb R)}\,.
\tag4.6
$$
We do this by setting
$$\aligned
\tilde T_0^{l,\omega}(f,g)&=
\sum_{\{j: \, 2^{\frac{3k}4}\le 2^j\le 2^{\frac{2k}3}\delta(l)^{-1}\}}
T^{l,\omega}(f,g_j)
\\
\\
\tilde T_1^{l,\omega}(f,g)&=\sum_{\{j: \, 2^j>2^{\frac{2k}3}\delta(l)^{-1}\}}
T^{l,\omega}(f,g_j)
\endaligned
$$
For the term $\tilde T^{l,\omega}_0(f,g)$, the index $j$ runs over
at most $|\log \delta(l)|$ terms. Thus, the
bound $(4.6)$ for this term results from the following
bound (uniform over $j$)
$$
\bigl\|\,T^{l,\omega}(f,g_j)\bigr\|_{L^2(dx\,dt)}
\le\,C\,\delta(l)^{-\frac 34}\,
|\log\delta(l)|^{\frac 12}\,
\,\|f\|_{L^2(\Bbb R^3)}\,\|g_j\|_{L^2(\Bbb R)}\,.
$$
This estimate follows from the argument for $(\frac 12,\frac 12)$ symbols
by freezing the $x_3$ coefficient in $a_j(t,x,\rho\omega)$, which
is possible now that the index $j$ is fixed.

To handle the term $\tilde T^{l,\omega}_1(f,g)$, we modify the argument
of Theorem 3.1 by taking the partition of unity such that
the symbol $a^{l,\omega}_{\nu,j,s,z}(t,x,\xi,\rho)$ is supported in the set
$$
\left|\,\frac{\xi}{|\xi|}-\xi^\nu\,\right|\le c\,2^{-\frac
k3}\,\delta(l)^{-1}
\,,\qquad
\bigl|\,\nabla_\xi\varphi_k(t,x,\omega)-z\,\bigr|\le 2\cdot2^{-\frac{2k}3}
\,,\qquad
|\,t-s\,|\le 2\cdot2^{-\frac{2k}3}\,\delta(l)^{-1}\,,
$$
and adjusting the spacing of the index points $(\nu,s,z)$ accordingly.
With these changes, and using the modified $S_{\frac 23,\frac 23}$
estimates for the symbol, estimates $(3.17)$ and $(3.19)$ are respectively
replaced by
$$
\matrix\format\r\;\;&\c&\;\;\l\\
\bigl\|\,
T^{l,\omega}_{\nu,j,s,z}\bigl(\,T^{l,\omega}_{\nu,j,s',z'}\,\bigr)^*
\,\bigr\|
&\le &\displaystyle
\frac{C_N\,\delta(l)^{-2}}
{\bigl(\,1+2^{\frac{2k}3}|\,z-z'|+
2^{\frac{2k}3}\,\delta(l)\,|\,s-s'|\,\bigr)^N}\\
\\
\\
\bigl\|
\,\bigl(T^{l,\omega}_{\nu,j,s,z}\,\bigr)^*\,
T^{l,\omega}_{\nu',j',s,z}\,\bigr\|
&\le&\displaystyle
\frac{C_N\,\delta(l)^{-1}}
{\bigl(\,1+2^{-\frac{2k}3}\,\delta(l)^2\,\bigl|\,2^j-2^{j'}\bigr|
+2^{\frac k3}\,\delta(l)\,|\,\xi^\nu-\xi^{\nu'}|\,
\bigr)^N}\,,
\endmatrix
$$
where we use the appropriate modification of Lemma 3.3.
Since the indices now run over $2^j\,\delta(l)\ge 2^{\frac{2k}3}$, the
rest of the proof of Theorem 3.1 goes through.
\qed\enddemo

In the case that $\Omega$ is the complement in $\R^3$ of a strictly
convex obstacle, with the Euclidean metric understood, a partition
of unity argument allows one to extend Theorem 4.1 to hold globally
on $\Omega$ (but still over a finite time interval.) Precisely, from the
result of Klainerman-Machedon \cite{2}
that the conclusion of the theorem holds globally on 
Minkowski space, together with finite propagation velocity and 
energy estimates, we may conclude the following extension.

\proclaim{Theorem 4.2}
Let $\Omega$ be the complement in $\R^3$ of a strictly convex,
smoothly bounded compact subset.
Suppose that $u$ and $v$ satisfy the Cauchy problem for the
Euclidean metric on $\Omega$
with Dirichlet condition, with respective data
$$
u_0,v_0\in H^2_D(\Omega)\,,\quad
u_1,v_1\in H^1_D(\Omega)\,,\quad
F,G,DF,DG\in L^1_t(\,[-1,1];L^2(\Omega))\,.
$$
Then the following hold, for any of the null forms $Q$,
$$
\multline
\bigl\|DQ(du,dv)\bigr\|_{L^2_{t,x}(\,[-1,1]\times\Omega)}\le
C\,
\Bigl(\,\|u_0\|_{H^2_D(\Omega)}+\|u_1\|_{H^1_D(\Omega)}
+\sum_{|\alpha|\le 1}
\|D^\alpha F\|_{L^1_tL^2_x(\,[-1,1]\times\Omega)}\,\Bigr)\\
\times
\Bigl(\,\|v_0\|_{H^2_D(\Omega)}+\|v_1\|_{H^1_D(\Omega)}
+\sum_{|\alpha|\le 1}
\|D^\alpha G\|_{L^1_tL^2_x(\,[-1,1]\times\Omega)}\,\Bigr)
\,.
\endmultline
$$
$$
\multline
\bigl\|Q(du,dv)\bigr\|_{L^2_{t,x}(\,[-1,1]\times\Omega)}\le
C\,
\Bigl(\,\|u_0\|_{H^1_D(\Omega)}+\|u_1\|_{L^2(\Omega)}
+\|F\|_{L^1_tL^2_x(\,[-1,1]\times\Omega)}\,\Bigr)\\
\times
\Bigl(\,\|v_0\|_{H^2_D(\Omega)}+\|v_1\|_{H^1_D(\Omega)}
+\sum_{|\alpha|\le 1}
\|D^\alpha G\|_{L^1_tL^2_x(\,[-1,1]\times\Omega)}\,\Bigr)
\,.
\endmultline
$$
\endproclaim

\enddocument